\documentclass[11pt]{amsart}
\usepackage{amssymb,amsfonts,amscd,verbatim}
\usepackage[leqno]{amsmath}
\usepackage{mathrsfs}
\usepackage{xy} 

\newcommand{\coker}{\operatorname{coker}} 
\newcommand{\Shv}{{\rm Shv}}
\newcommand{\Sh}{{\rm Sh}}
\newcommand{\Ext}{{\rm Ext}}
\newcommand{\Hom}{{\rm Hom}}
\newcommand{\Spec}{{\rm Spec}}
\newcommand{\by}[1]{\stackrel{#1}{\rightarrow}}

\newcommand{\G}{\mathbb{G}}  
\newcommand{\C}{\mathbb{C}}  
\newcommand{\Z}{\mathbb{Z}}  

\newcommand{\df}{\mbox{\,${:=}$}\,}

\newcommand{\ie}{{\it i.e.},\ }
\newcommand{\cf}{{\it cf.\/}\ }

\newcommand{\op}{{\it op. cit.\/} }

\newcommand{\id}{{\rm id}}
\newcommand{\eff}{{\rm eff}}
\newcommand{\et}{{ \rm {\acute et}}}
\newcommand{\tor}{{\rm tor}}

\newcommand{\fr}{{\rm fr}}

\newcommand{\fppf}{{\rm fppf}}
\newcommand{\M}{\mathcal{M}_1} 
\newcommand{\tM}{{}^t\!\mathcal{M}_1} 

 \swapnumbers

\numberwithin{equation}{section}
\newtheorem{thm}[equation]{Theorem} 
\newtheorem{propose}[equation]{Proposition}
\newtheorem{lemma}[equation]{Lemma}

\theoremstyle{definition}
\newtheorem{defn}[equation]{Definition}
\newtheorem{remark}[equation]{Remark}
\newtheorem{example}[equation]{Example}

\begin{document}

\input xy  
\xyoption{all}

\title{Remarks on $1$--motivic sheaves}
\author{Alessandra Bertapelle}
\address{Dipartimento di Matematica, Universit\`a degli Studi di
Padova\\Via Trieste, 63\\Padova I-35121\\ Italy}
\email{alessandra.bertapelle@unipd.it}
\begin{abstract} {\normalsize We generalize the construction of the category of $1$--motives with torsion $\tM$ in \cite{BVRS} as well as the construction of the
category of $1$--motivic sheaves $\Shv_1$ in
\cite{BVK} to perfect fields $k$ (without inverting the exponential
characteristic). For $k$ transcendental over the prime field we extend a result
in \cite{BVK} showing that $\tM$ and $\Shv_1$ have equivalent bounded derived
categories.} \end{abstract}

\date{May 24, 2012}

\maketitle
\section*{Introduction}
Let $k$ be a field of characteristic $0$. 
In \cite{BVK} the authors introduce the category of $1$--motivic sheaves $\Shv_1$ and show that $D^b(\Shv_1)$ and the bounded derived category of $1$--motives with torsion $D^b(\tM)$ are both equivalent to the thick subcategory of Voevodsky's triangulated category of motives ${\rm DM}^\eff_{gm}(k)$ generated by motives of smooth curves. 
When $k$ is a perfect field of positive characteristic $p$, the definition of $\tM$ in \cite{BVRS} does not work well: for example, it fails to provide an abelian category.
However, it is sufficient to invert $p$-multiplications to get a $\Z[1/p]$-linear abelian category (\cite{BVK}, C.5.3). 
On the other hand, the general definition of the category of $1$--motivic sheaves in \cite{BVK} also requires inverting $p$-multiplications and an equivalence result between bounded derived categories still holds (\cite{BVK}, 3.9).
As the authors explain, if one is interested in comparison results with Voevodsky's category {\rm DM}, one can not avoid inverting the exponential characteristic.

In this paper we show that an integral definition of $\tM$ is possible over any perfect field $k$, \ie without inverting the exponential characteristic, if we allow finite connected $k$-group schemes in the component of degree $-1$ of $1$--motives. 
We then get an abelian category, still denoted by $\tM$, that contains the category of Deligne $1$--motives as a full exact subcategory.
Moreover, on passing to the fppf topology, the definition of $1$--motivic sheaves works ``integrally'' over any perfect field providing the abelian category of $1$--motivic (fppf) sheaves $\Shv_1^\fppf$. For $k$ of characteristic $0$ both constructions are equivalent to those given in \cite{BVK}. Furthermore, for $k$ perfect of positive characteristic and transcendental over its prime field, we show that the categories $D^b(\Shv_1^\fppf)$ and $D^b(\tM)$ are equivalent, \ie we get an integral version of \cite{BVK}, 3.9.2, without passing through Voevodsky's category DM (\cf~\ref{thm.main} \& \ref{pro.voe}). 
Our proof only works for fields of characteristic $0$ and perfect fields of positive characteristic that are transcendental over their prime fields because it requires the existence of semiabelian varieties having non-torsion points. It is not clear if the comparison result in \cite{BVK}, 3.9.2, holds over a finite field without inverting the exponential characteristic.

In the final part of the paper we give a useful characterization of $1$--motivic sheaves as quotients of $k$-group schemes. 

This paper is a shortened version of the earlier preprint posted to arXiv.
Analogous results for Laumon $1$-motives will be studied in \cite{BER}.

\section{The derived category of $1$--motives with torsion}
In this section we generalize the construction of the category of $1$--motives with torsion
$\tM$ (\cf~\cite{BVRS}, \cite{BVK}) and we study its subcategory $\M^\star$.

Let $k$ be a perfect field and $p$ its exponential
characteristic. We say that a $k$-group scheme is \emph{discrete} if
it is finitely generated locally constant for the \'etale topology
(\cf~\cite{BVK}, 1.1.1). Let $\mathcal{CE}$ be the category of
commutative $k$-group schemes $\Gamma$ that are products of a
discrete group scheme $\Gamma^\et$ by a finite commutative connected
group scheme $\Gamma^0$. The group $\Gamma$ contains a maximal
torsion subgroup $\Gamma_\tor$, so that
$\Gamma_\fr:=\Gamma/\Gamma_\tor$ is discrete and torsion free.
Observe that it follows from the perfectness of $k$ that
$\mathcal{CE}$ is closed under extensions. It suffices to consider
the case of an extension $\Gamma$ of a discrete group scheme
$\Gamma^\et$ by a finite commutative connected group scheme
$\Gamma^0$ and the case of an extension of $\Gamma^0$ by
$\Gamma^\et$. For $\Gamma^\et$ finite, the result follows from the
classification of finite commutative $k$-group schemes; in the case
$\Gamma^\et$ torsion free, if $n$ is the order of $\Gamma^0$, 
there are surjective homomorphisms
$0=\Ext(\Gamma^\et/n\Gamma^\et,\Gamma^0)\to
\Ext(\Gamma^\et,\Gamma^0)$ and
$0=\Hom(\Gamma^0,\Gamma^\et/n\Gamma^\et)\to
\Ext(\Gamma^0,\Gamma^\et)$.

\begin{defn}\label{def.effm}
An \emph{effective $1$--motive with torsion} is a complex of
$k$-group schemes $M=[u\colon L\to G]$ where $L$ is an object in
$\mathcal{CE}$ and $G$ is semiabelian. An \emph{effective morphism}
$M\to M'$ is a map of complexes $(f,g)$, with $f\colon L\to L'$,
$g\colon G\to G'$ morphisms of $k$-group schemes. $M$ is said to be
\emph{\'etale} if $L$ is \'etale.
\end{defn}
Denote by
$\tM^{\eff}$ the category of effective $1$--motives with torsion and by
$\tM^{\eff,\et}$ the full subcategory of \'etale effective
$1$--motives with torsion. The category $\M$ of Deligne $1$--motives is the full
subcategory of $\tM^{\eff,\et}$ consisting of those $M=[L\to G]$ with $L$
torsion-free (\cf~\cite{DEL}).

\begin{defn} Let $\Sigma$ be the class of quasi-isomorphisms (q.i) of effective
 $1$--motives with torsion, \ie the class of effective maps
$(f,g)\colon M\to M'$ where $g$ is an isogeny, $f$ is surjective
and $\ker(f)=\ker(g)$ is a finite group scheme. The category of
\emph{$1$--motives with torsion} $\tM$ is then defined as
the localization of $\tM^{\eff}$ at $\Sigma$. Similarly,
the category of \emph{\'etale} $1$--motives with torsion
$\tM^{\et}$ is the localization of $\tM^{\eff,\et}$ at the class of quasi-isomorphisms.
\end{defn}

For $k$ of characteristic $0$ the category $\tM^\et$ was first introduced in 
\cite{BVRS} where it was denoted by $\M$,
 and then over any perfect field in
\cite{BVK} where it was denoted by $\tM$. It was
proved to be equivalent to the category of Mixed Hodge Structures
of level $\leq 1$ for $k=\C$ (\cf~\cite{BVRS}, 1.5). For $k$ of
characteristic zero, $\tM^\et=\tM $ is an abelian category 
(\cf~\cite{BVRS}, \cite{BVK}). Over a field $k$ of positive
characteristic $p$, the category $\tM^\et$
 is a subcategory of $\tM$ that becomes
abelian upon inversion of $p$-multiplications (\cf~\cite{BVK}, C.5.3).
Observe that the natural (faithful) functor $\tM^\et\to \tM$
 has a left inverse that associates to a $1$--motive with torsion
$M=[u\colon \Gamma\to G]$ the \'etale $1$--motive with torsion
$M^\et=[\Gamma/\Gamma^0\to G/u(\Gamma^0)]$. The effective map $M\to
M^\et$ becomes an isomorphism upon inversion of $p$-multiplications
and hence the abelian category $\tM^\et[1/p]$ is equivalent to
$\tM[1/p]$. In the sequel we will prove that $\tM$ itself is an
abelian category.

\begin{remark}
To motivate Definition~\ref{def.effm}, observe that in
the characteristic zero case,
starting with an effective $1$--motive with torsion $M=[L\to G]$ and an isogeny
$g\colon G'\to G$, by pull-back one always gets a q.i.
$(f,g)\colon [L'\to G']\to M $.
Over fields of positive characteristic this is not always the case if $L$
is required to be a discrete group because there are isogenies with connected
kernel. Therefore one can generalize the construction of the category of
$1$--motives with torsion in \cite{BVRS} either by inverting
$p$-multiplications, as done in \cite{BVK},
or by allowing non-\'etale finite group schemes, as we do here.
\end{remark}
 To show that $\tM$ is an abelian category, we will
follow the analogous proof for $\tM^\et[1/p]$ in \cite{BVK},
Appendix~C. One can check that all results in Appendix C.1-C.5
of \cite{BVK} extend to our more general context. For the reader's convenience
we present some proofs.

\begin{lemma}\label{lem.limcolim}
The category of Deligne $1$--motives $\M$ and
the category of effective $1$--motives with torsion $\tM^{\eff}$
both admit finite limits and colimits.
\end{lemma} \proof (\cf~\cite{BVK}, C.1.3.)
Since we are working with additive categories, it suffices to show
that our categories have kernels and cokernels. For the definition
of the kernel of an effective morphism $\varphi=(f,g)$ take
$\ker(\varphi)=[\ker^0(f)\to \ker^0(g)]$ where $\ker^0(g)$ is the
reduced subgroup of the identity component of the kernel (as group
schemes) of $g$ and $\ker^0(f)$ is the pull-back of $\ker^0(g)$
along $\ker(f)\to \ker(g)$; the definition makes sense because
$k$ is perfect. The cokernel of $\varphi$ is the cokernel as group
schemes in each degree, \ie $[\bar u\colon \coker (f)\to \coker(g)]$.
Clearly, if $\varphi$ is a morphism of Deligne
$1$--motives, $\ker^0(f)$ is \'etale and torsion free. As $\coker(f)$ is \'etale
but in general not torsion free one takes as cokernel of $\varphi$ in $\M$ the
$1$--motive $[\coker(f)_\fr\to \coker(g)/\bar u(\coker(f)_\tor)]$.
 \qed

We will see now that the class of q.i.\ admits a calculus of right
fractions; hence any morphism $\varphi\colon M\to M'$ in the
category $\tM$ can be represented by a (not unique) pair
$(\sigma,\psi)$ where $\sigma\colon M''\to M$ is a q.i., $\psi\colon
M''\to M'$ is an effective map and $\varphi\sigma=\psi$ in $\tM$.

\begin{lemma}\label{lem.simpl} Morphisms in $\Sigma$ are
simplifiable both on the left and on the right.
\end{lemma}
\proof (\cf~\cite{BVK}, C.2.3.)
Let $(f,g)\in \Sigma$. Since $f\colon L\to L'$, $g\colon G\to G'$ are
epimorphisms, they are simplifiable on the
right. Suppose given an effective map $(f',g')$ with $ff'=0$, $gg'=0$. Since $g$ is an isogeny,
say of degree $n$, the multiplication by $n$ on $G$ factors through $g$ and
$ng'=0$. Hence $g'=0$; furthermore $f'=0$ because $f'$, $g'$ both factor
through $\ker(g)=\ker(f)$.
\qed

\begin{lemma}\label{lem.calculus} $\Sigma$ is a left multiplicative system and
admits a calculus of right fractions.
\end{lemma}
\proof (\cf~\cite{BVK}, C.2.4.) Let $\varphi'\colon M''\to M'$ be a
q.i.\ and $\varphi\colon M\to M'$ an effective map. For the first
assertion, we have to show that there exist an effective map $\psi\colon \tilde
M \to M''$ and a q.i.\ $\psi'\colon \tilde M\to M$ forming a commutative diagram
with $\varphi$ and $\varphi'$. Define $\tilde G$ as the reduced
subgroup of the identity component of $G\times_{G'} G''$; since $k$ is
perfect, $\tilde G$ is still an algebraic $k$-group and indeed it is a semiabelian variety
isogenous to $G$. Define now
 $\tilde M\df [\tilde L\to \tilde G]$ via pull-back and let $\psi,\psi'$ be the
obvious maps. Since $\Sigma$ is a left multiplicative system and
its maps are simplifiable on the left by Lemma~\ref{lem.simpl}, it admits a calculus of right
fractions. \qed

In order to prove that $\tM$ is an abelian category, the following
notion of strict morphism will play a key role.

\begin{defn}
An effective morphism $\varphi=(f,g)\colon M\to M'$ is \emph{strict}
if $g$ has smooth connected
 kernel, \ie the kernel of $g$ is still semiabelian.
\end{defn}

In particular if $\varphi$ is strict and $g$ is an isogeny then $g$ is an isomorphism.

\begin{propose}\label{pro.fact}
Any effective morphism $\varphi=(f,g)\colon M\to M'$ factors as
$\sigma\tilde\varphi=\varphi$ with $\sigma$ a quasi-isomorphism and
$\tilde\varphi\colon M\to \tilde M$ strict.
\end{propose}
\proof (\cf~\cite{BVRS}, 1.3.)
We recall the main
ideas of the proof. Let $M=[L\to G]$ and $M'=[L'\to G']$.
It is sufficient to work with the morphism of semiabelian
varieties $g\colon G\to G'$.
Up to dividing $G$ by $\ker^0(g)$, we may assume that $\ker(g)$
 is finite of order $n$. Define $\tilde G:=G'/g({}_n G)$ with
${}_n G$ the kernel of the $n$-multiplication on $G$. Writing down a
diagram with horizontal sequences the short exact sequences of
$n$-multiplication for $G$ and $G'$, and vertical sequences the exact
sequence \[0\to \ker(g)\to G\to G'\to 0\] twice, one checks that the
$n$-multiplication of $G'$ factors through an isogeny
$h\colon\tilde G\to G'$. Moreover $g\colon G\to G'$ lifts to a
monomorphism $\tilde g\colon G\to \tilde G$ so that $g=h\tilde
g$. The $1$-motive $\tilde M$ is then defined via the pull-back of $L'$ along
$h$.
\qed

\begin{example}\label{ex.gm}
Let $n$ be a positive integer, consider the $n$-multiplication on $\G_m$ and let ${\boldsymbol \mu}_n$ denote its kernel. The above maps factors as
\[\G_m\by{(0,\id)} [{\boldsymbol \mu}_n\to \G_m]\by{(0,n)} \G_m\]
where the first map is a strict morphism and the second one is a
quasi-isomorphism. Observe that the $1$--motive in the middle is not
a $1$--motive with torsion in the sense of \cite{BVRS} if $n$ is
not invertible in $k$.

Furthermore, the cokernel in $\tM^\eff$ of the $1$--motive $n\colon \G_m\to
\G_m$ is trivial, while the cokernel of $(0,\id)\colon \G_m\to [{\boldsymbol
 \mu}_n\to \G_m]$
is $[{\boldsymbol \mu}_n\to 0]$. Since the map \[(0,n)\colon [{\boldsymbol
 \mu}_n\to \G_m]\to \G_m\]
 is an isomorphism in
$\tM$, it is clear that we can not expect the functor $\M^\eff\to \tM$ to
preserve cokernels. However it does preserve cokernels for strict morphisms.
\end{example}

\begin{lemma}\label{lem.coker}
Let $\varphi\colon M\to M'$ be a strict effective morphism. Its
cokernel $\coker(\varphi)$ in the category $\tM^\eff$ remains a
cokernel in $\tM$.
\end{lemma}
\proof (\cf~\cite{BVK}, C.5.2.) Let
 $\varphi\colon M=[L\to G]\to M'=[L'\to G']$ and let
 \mbox{$\psi\colon M'\to M''$} be a
morphism of $1$--motives with torsion such that $\psi\varphi=0$ in $\tM$. We
check that $\psi$
factors uniquely through $\coker(\varphi)=[L'/L\to G'/G]$, the cokernel of $\varphi$ in
the category $\tM^\eff$. Let $s\colon N'\to M'$
be a q.i.\ so that $\psi s$ is effective. Let $\varphi'\colon N\to
N'$ be an effective map and $t\colon N\to M$ a q.i.\ such that
$s\varphi'=\varphi t$. Note that they exist by the
calculus of right fractions. It
is not difficult to check that $s,t$ induce a q.i.\ $\coker
(\varphi')\to \coker(\varphi)$. Since $\psi s \varphi'=0$, the
effective map $\psi s$ factors uniquely through $\coker (\varphi')$
and hence $\psi$ too factors through a morphism $\coker(\varphi)\to
M''$ in $\tM$. The uniqueness of factorization is proved with
similar arguments. \qed

By contrast, the functor
$\tM^{\eff}\to \tM$ always preserves kernels.

\begin{lemma}\label{lem.leftexact}
Kernels exist in $\tM$. The canonical functor $\tM^{\eff}\to \tM$ is
left exact and faithful.
\end{lemma}
\proof (\cf~\cite{BVK} C.5.1.)
The faithfulness follows form
Lemma~\ref{lem.simpl}. The existence of kernels in $\tM$ and the left
exactness 
of the above functor follow
 from the existence of kernels in $\tM^\eff$ (\cf~Lemma~\ref{lem.limcolim})
 and Lemma~\ref{lem.calculus}.
\qed

We now prove the main theorem on $1$--motives with torsion:

\begin{thm}\label{thm.ac}
\begin{itemize}
\item[(i)] The category $\tM$ is abelian.
\item[(ii)] Every short exact sequence of $1$--motives in $\tM$
\[ 0\to M'\to M\to M''\to 0
\]
can be represented (up to isomorphisms) by a sequence of
effective $1$--motives that is exact as sequence of complexes.
\item[(iii)] The natural functor $\M\to \tM$ from the category of Deligne
 $1$--motives to the category of $1$--motives with torsion is fully faithful and makes
 $\M$ an exact subcategory of $\tM$. In particular,
given a short exact sequence as in (ii) with
$M', M''$ isomorphic to Deligne $1$--motives, $M$ is
isomorphic to a Deligne $1$--motive (unique up to isomorphisms of Deligne $1$--motives).
\end{itemize}
\end{thm}
\proof (\cf~\cite{BVK}, C.5.3, C.7.1 \& \cite{BVRS}, 1.3.)
We start by proving (i). Kernels exist by
Lemma~\ref{lem.leftexact}. For the existence
of cokernels, by Proposition~\ref{pro.fact}, we may work with strict
effective morphisms $\varphi\colon M\to M'$. By
Lemma~\ref{lem.coker} the cokernel of $\varphi$ in $\tM$ exists. Since
\[\coker(\ker\varphi\to M)\to \ker(M'\to \coker(\varphi))\] is an
isomorphism in $\tM^\eff$, it is the same in $\tM$ and we are done.
Statement (ii) follows immediately from the proof of (i) on taking a
strict ``representative'' of $M'\to M$. For (iii): Clearly the
functor $\M\to \tM$ is fully faithful because it has
 a left inverse/left adjoint sending
$M=[u\colon L\to G]$ to $M_\fr:=[L_\fr\to G/u(L_\tor)]$. By (ii) any
short exact sequence of $1$--motives in $\tM$ is represented up to
isomorphisms by a sequence of effective $1$--motives that is exact
as sequence of complexes. Since the outer $1$--motives are
isomorphic to Deligne $1$--motives, by a direct computation one sees
that the same holds for
 the one in the middle. \qed

All information needed to understand the bounded
derived category of $\tM$ is encoded in $\M$ and indeed in
the following subcategory $\M^\star$, if $k$ is large enough.

\begin{defn}
Denote by $\M^\star$ the full subcategory of the category of Deligne $1$--motives $\M$
 whose objects are effective $1$--motives
$[u\colon L\to G]$ with $\ker u=0$ (and $L$ \'etale torsion free).
\end{defn}

 \begin{remark}\label{rem.ker}
 Observe that there are no non-trivial q.i.\ in $\M^\star$ and that,
 by Theorem~\ref{thm.ac} (iii), the natural functor $\M^\star\to \tM$ is fully
 faithful; hence we may treat $\M^\star$ as a full subcategory of $\tM$.
 Moreover, given two quasi-isomorphic $1$--motives with torsion $M_i=[u_i\colon L_i\to G_i]$,
$i=1,2$, $\ker(u_1)$ is trivial if and only if $\ker(u_2)$ is
trivial. In particular, $M=[u\colon L\to G]$ is quasi-isomorphic to
a $1$--motive in $\M^\star$ if and only if $\ker(u)=0$.

If $k$ is a finite field, or more generally a field of positive
characteristic and algebraic over its prime field, then $\M^\star$
is equivalent to the category of semiabelian varieties. Indeed there
are no injective maps from $\Z$ into a semiabelian variety. Under this
hypothesis, the category $\M^\star$ will be of no use for describing
$D^b(\tM)$.
\end{remark}

\begin{lemma}\label{lem.mw}
The category of Deligne $1$--motives $\M$ is a generating subcategory of $\tM$ closed under kernels
and closed under extensions.
Suppose that $k$ has characteristic $0$ or is transcendental over
its prime field. Then $\M^\star$ too is a generating subcategory
of $\tM$ closed under both kernels and extensions.
\end{lemma}
\proof By the description of kernels in
Proposition~\ref{lem.limcolim} and Lemma~\ref{lem.leftexact}, both
$\M^\star$ and $\M$ are closed under kernels.

Now $\M$ is closed under extensions by Theorem~\ref{thm.ac}~(ii).
For $\M^\star$: let
 \[M^\bullet=0\to M'\to M\to M''\to 0\]
 be a short exact sequence in $\tM$
 with $M',M''$ in $\M^\star$, by Theorem~\ref{thm.ac}~(ii) and Remark~\ref{rem.ker}
 the effective $1$--motive $M=[u\colon L\to G]$ is, up to quasi-isomorphisms, a Deligne
$1$--motive $\tilde M=[\tilde u\colon \tilde L\to \tilde G]$
and $\ker(u)=0$. Hence also $\ker(\tilde u)=0$ and $\tilde M$ is an object of $\M^\star$.

To see that $\M$ (respectively $\M^\star$) is generating, we have to show
that for any $1$--motive with torsion $M=[u\colon L\to G]$ there
exists an epimorphism $\varphi\colon M'\to M$ with $M'$ in $\M$
(respectively \ in $\M^\star$). We may assume that $L$ is \'etale. Indeed,
if $L=L^0\times L^\et$ we may always find an embedding $v\colon
L^0\to B$ of $L^0$ into an abelian variety $B$. Then the $1$--motive
$N=[w\colon L^0\times L^\et \to B\times G]$, with
$w(x,y)=(v(x),u(x,y))$, is quasi isomorphic to the \'etale effective
$1$--motive $N'=[L^\et\to B\times G/w(L^0)]$ because $L^0$ embeds
into $B\times G$. Furthermore the effective morphism $N\to M$ that
is the identity in degree $-1$ and the projection $B\times G\to G$
in degree $0$ is a strict epimorphism, hence it remains an
epimorphism in $\tM$ (\cf~Lemma~\ref{lem.coker}). 
In particular, we have an epimorphism $N'\to M$
in $\tM$.

Assume then $L$ \'etale. If it is torsion free, $M$ is already a
Deligne $1$--motive. Otherwise
 choose an epimorphism $g\colon L'\to L$, with $L'$ discrete torsion free,
and consider the Deligne $1$--motive $M'=[u\circ g\colon L'\to G]$
and the strict epimorphism $(g,\id)\colon M'\to M$. The latter remains an
epimorphism in $\tM$, which concludes the proof that $\M$ is a
generating subcategory of $\tM$.

To see that $\M^\star$ is generating, there remains to verify that there exists an
epimorphism $\varphi\colon M''\to M'$ with $M''$ in $\M^\star$ and
$M'=[u\circ g\colon L'\to G]$
the $1$--motive constructed above.
Observe that there
exists a torus $H$ and a monomorphism $h\colon L'\to H$: one may work over
 a finite separable extension $k'$ of $k$ so that the base change of $L'$ to
 $k'$ is isomorphic to $\Z^d$. A
 monomorphism $ L'\times_k {k'}\to \G_m^d$ exists over $k'$
 because of the hypothesis on $k$. Then one descends this monomorphism to $k$
 by restriction of scalars.
Consider now the $1$--motive $M''\df [u''\colon L'\to H\times G]$,
$u''=(h,u)$, that is an object of $\M^\star$, and the morphism
$\varphi=(\id,p_G)\colon M''\to M'$ with $p_G$ the usual projection
map. Since $\varphi$ is a strict epimorphism, it remains an epimorphism
in $\tM$ by Lemma~\ref{lem.coker} and we are done. \qed

\section{The category of $1$--motivic sheaves}\label{sec.msh}

In this section we introduce the category of $1$--motivic sheaves $\Shv_1^\fppf$
 and we study its subcategory $\Shv_1^\star$.

\begin{defn}\label{def.motsh}
A sheaf $\mathcal F$ for the fppf topology
over $\Spec(k)$ is \emph{$1$--motivic} if there exists a morphism
of sheaves $b\colon G\to \mathcal F$ with $G$ a semiabelian variety
over $k$ and $\ker (b)$, $\coker (b)$ in $\mathcal{CE}$.
 The map
$b$ is said to be \emph{normalized} if $\ker (b)$ is \'etale and
torsion-free.

Denote by $\Shv_1^\fppf$ the full subcategory of the category of
fppf sheaves over $\Spec(k)$ whose objects are the $1$--motivic
sheaves. Denote by $\Shv_1^\fr$ the full subcategory of those
$1$--motivic sheaves with both $\ker(b)$ and $\coker(b)$ discrete
torsion free.
\end{defn}

Observe that we have an exact sequence
\begin{eqnarray}\label{eq.F}
0\to L\by{a} G \by{b}\mathcal F\by{c} E\to 0
\end{eqnarray}
with $L$ and $E$ in $\mathcal{CE}$. For $k$ of characteristic $0$,
the category $\Shv_1^\fppf$ is equivalent to the category $\Shv_1$
defined in \cite{BVK} (\cf~\op 3.3.2). We will explain in
Section~\ref{sec.result} the relation between $\Shv_1$ and our
$\Shv_1^\fppf$ over general perfect fields. Denote by $\Shv_0^\fppf$
the full subcategory of $\Shv_1^\fppf$ consisting of those $\mathcal
F$ with $G=0$; it is equivalent to $\mathcal{CE}$.

\begin{propose}\label{pro.motsh}
\begin{itemize}
\item[(i)] In Definition~\ref{def.motsh} we may choose $b$ normalized.
\item[(ii)] Given two $1$--motivic sheaves $\mathcal F_1, \mathcal F_2$,
normalized morphisms $b_i\colon G_i\to \mathcal F_i$, $i=1,2$, and
a morphism of sheaves
 $\varphi\colon \mathcal F_1\to \mathcal F_2$
 there exists a unique homomorphism of group schemes
$\varphi_G\colon G_1\to G_2$ such that $\varphi\circ b_1=b_2\circ \varphi_G$.
\item[(iii)] Given a $1$--motivic sheaf $\mathcal F$, a
 morphism $b\colon G\to \mathcal F$ as above with $b$ normalized
is uniquely (up to isomorphisms) determined by $\mathcal F$.
\item[(iv)]
$\Shv_1^\fppf$ and $\Shv_0^\fppf$ are exact abelian subcategories
of the category of fppf sheaves over $\Spec(k)$.
\end{itemize}
\end{propose}\proof
 (\cf~\cite{BVK}, 3.2.3). (i). If $b$ is not normalized one
simply divides $G$ by $L_\tor$ for $L=\ker(b)$.

For (ii) consider the sequence
\eqref{eq.F} for both $\mathcal F_i$, $i=1,2$. As $c_2\varphi
b_1=0$, the morphism $\varphi$ induces a morphism $\varphi_E\colon
E_1\to E_2$ and hence a morphism $G_1\to G_2/L_2$. Since
$\Ext^1(G_1,L_2)=0$ (\cf~\cite{GRO2}, VIII, 3.4.2 \& 5.5),
we get a morphism \mbox{$\varphi_G\colon G_1\to G_2$} and hence a
\mbox{$\varphi_L\colon L_1\to L_2$}. The map $\varphi_G$ is unique because
$\Hom(G_1,L_2)=0$. (iii) follows from (ii). The assertion in (iv) is
clear for $\Shv_0^\fppf$. For $\Shv_1^\fppf$ one checks directly
that the kernel and the cokernel of a morphism of $1$--motivic
sheaves \mbox{$\varphi\colon {\mathcal F}_1\to {\mathcal F}_2$} are still
$1$--motivic: one considers $b_3\colon G_3\to \ker(\varphi)$ with
$G_3$ the reduced subgroup of the identity component of the kernel
of $\varphi_G\colon G_1\to G_2$,
 and \mbox{$b_4\colon \coker(\varphi_G)\to \coker(\varphi)$}. An easy diagram chase
shows that kernels and cokernels of $b_3$ and $b_4$ are in $\mathcal{CE}$. \qed

We now introduce a full subcategory of the category of $1$--motivic
sheaves that will play the ``dual role'' of $\M^\star$ for $\tM$.

\begin{defn} Denote by $\Shv_1^\star$ the full subcategory of $\Shv_1^\fppf$
consisting of those sheaves $\mathcal F$ such that there exists a
$b\colon G\to \mathcal F$ with
 $b$ an epimorphism, \ie $E=\coker(b)=0$.
\end{defn}

For $k$ algebraic over a finite field, $\Shv_1^\star$ is equivalent to the
category of semiabelian varieties while $\Shv_1^\fppf$ is the category of those fppf
sheaves that are extensions of a sheaf $E$ in $\mathcal{CE}$ by a semiabelian variety; in particular
they are always representable by a $k$-group scheme (\cf~\cite{EtC}, III, 4.3).

 \begin{lemma}\label{lem.f0}
For a $1$--motivic sheaf $\mathcal F$ there exist unique (up to
isomorphisms)
 $\mathcal F^\star$ in $\Shv_1^\star$, $E$ in $\Shv_0$ and an exact sequence
 \[0\to \mathcal F^\star\to \mathcal F\to E\to 0.\]
 \end{lemma}
 \proof
 With notations as in \eqref{eq.F}, simply take $E=\coker (b)$ and $\mathcal F^\star=\coker(a)$.
 \qed

Consider now the following diagram:
\begin{eqnarray}\label{pic.catsh}
\xymatrix{
\M \ar[r]^d&\tM \ar[r]^{{\mathrm H}_0}& \Shv_1^\fppf \\
\M^\star \ar[u]^\iota \ar[rr]^{{\mathrm H}_0}& &\Shv_1^\star\ar[u]\\
 }
\end{eqnarray}
where the vertical arrows and $d$ are the usual inclusion functors.
The functor ${\mathrm H}_0$ maps a $1$--motive with torsion
$M=[u\colon L\to G]$ to the $1$--motivic sheaf $\coker (u)={\mathrm
H}_0(M)$. This sheaf is always in $\Shv_1^\star$. 

The following results motivates the introduction of the subcategory $\Shv_1^\star$.

\begin{lemma}\label{lem.a}
The functor ${\mathrm H}_0$ provides an equivalence of categories
between $\M^\star$ and $\Shv_1^\star$.
\end{lemma}
\proof It follows immediately from Proposition~\ref{pro.motsh} that
${\mathrm H}_0\colon \M^\star \to \Shv_1^\star$ is fully faithful.
It is also essentially surjective. Indeed, given a $1$--motivic sheaf
$\mathcal F$ in $\Shv_1^\star$ and a normalized morphism $b\colon
G\to \mathcal F$ as in \eqref{eq.F}, the Deligne $1$--motive
$[u\colon \ker (b)\to G]$ satisfies $\ker (u)=0$ and $\coker
(u)=\mathcal F$. \qed

The ``dual" of Lemma~\ref{lem.mw} holds:

\begin{lemma}\label{lem.shv10}
The category $\Shv_1^\fr$ is a cogenerating subcategory of $\Shv_1$
closed under cokernels and closed under extensions. If $k$ has
characteristic $0$ or if it is transcendental over its prime
field, the same holds for $\Shv_1^\star$.
\end{lemma}
\proof The only non-trivial fact is that the subcategory
$\Shv_1^\fr$ (respectively $\Shv_1^\star$) is cogenerating, \ie that for
any $1$--motivic sheaf $\mathcal F$ as in \eqref{eq.F} there exists
a $\mathcal F'$ in $\Shv_1^\fr$ (respectively in $\Shv_1^\star$) and a
monomorphism $\varphi\colon \mathcal F\to \mathcal F'$.

Suppose $b\colon G\to \mathcal F$ in \eqref{eq.F} normalized. 
Let $\tilde{ \mathcal F}$ be
the pull-back of $\mathcal F$ along $E_\tor\to E$. We have $\tilde
b\colon G\to \tilde {\mathcal F}$ whose kernel is still $L$ and the
cokernel is $E_\tor$. Furthermore we have a short exact sequence
\[0\to \tilde{\mathcal F}\to \mathcal F\to E_\fr\to 0.\]
To show that $\Shv_1^\fr$ is cogenerating, it is sufficient to show
that $\tilde{\mathcal F}$ embeds into a sheaf $\mathcal G$ in
$\Shv_1^\star$: the push-out of $\mathcal F$ along $\tilde {\mathcal F}
\to \mathcal G$ provides then a monomorphism $\mathcal F\to
{\mathcal F}'$ with $\mathcal F'$ in $\Shv_1^\fr$ because extension
of $E_\fr$ by $\mathcal G$.

By Proposition~\ref{pro.van}, $\tilde {\mathcal F}$ is $\coker(L\to
F_0)$ with $F_0$ an extension of $E_\tor$ by $G$. We show that
$F_0$ embeds into a semiabelian variety $\tilde G$: Denote by ${}_n
F_0$ the kernel of the $n$-multiplication on $F_0$, and similar
notation for ${}_nG$. Then ${}_n F_0$ is extension of $E_\tor$ by
${}_nG$, hence finite, and $F_0/{}_n F_0\cong G$. Choose an
embedding $f\colon {}_n F_0\to B$ into an abelian variety $B$ and
let $f'\colon F_0\to \tilde G$ be the push-out of $F_0$ along $f$.
Since $\tilde G$ is an extension of $G$ by $B$ it is a semiabelian
variety. Define now $\mathcal G=\coker(L\to \tilde G)$ and choose
the embedding $\tilde \varphi\colon \tilde{ \mathcal F}\to \mathcal
G$ induced by $F_0\to \tilde G$. As explained above $\tilde \varphi$
provides a monomorphism $\mathcal F\to {\mathcal F}'$ with
${\mathcal F}'$ in $\Shv_1^\fr$.

To show that $\Shv_1^\star$ is cogenerating we may suppose then that
$\mathcal F$ is in $\Shv_1^\fr$, \ie $E$ is \'etale and torsion free.
Consider the diagram \eqref{dia.def} describing $\mathcal F$ as
cokernel of a morphism $F_1\to F_0$ where $F_0$ is extension of a
discrete group $\Gamma=F_0/G$ by $G$. Such an extension splits over a
suitable finite separable extension $k'$ of $k$. Let $f\colon
\Spec(k')\to \Spec(k)$ be such that $f^*F_0=F_{0,k'}$ is isomorphic
to $G_{k'}\times \Gamma_{k'}$. We may assume that $\Gamma$ is
constant over $k'$. Embed $\Gamma_{k'}$ into a semiabelian 
variety $H_{k'}$ over $k'$. This is possible because of the
assumptions on $k$. We get then a monomorphism $\Gamma\to
f_*H_{k'}$, where $f_*H_{k'}$ is still a semiabelian variety, namely
the Weil restriction of $H_{k'}$. Moreover, $G\to f_*f^*G$ is a
monomorphism. Hence we have a composition of monomorphisms $F_0\to
f_*f^* F_0\to (f_*f^* G)\times f_*H_{k'}$ where the latter scheme is
semiabelian. Define now $\mathcal F'$ in $\Shv_1^\star$ as the
cokernel of $F_1\to (f_*f^* G)\times f_*H_{k'}$; by construction
$\mathcal F$ embeds into $\mathcal F'$. \qed

\section{Equivalence on bounded derived categories}

In this section we derive results on the
 bounded derived categories of $\tM$ and $\Shv_1^\fppf$.
For conciseness let $(*)$ denote the following hypothesis:
\begin{equation*}(*)\quad k\ \text{\it has characteristic zero or it is transcendental over
 its prime field.}
\end{equation*}
 In order to see that ${\mathrm H}_0$ induces an equivalence between the bounded
derived categories of $\tM$ and $\Shv^\fppf$, we check the following
facts:

\begin{lemma}\label{lem.der1}
Assume condition $(*)$. 
Let $N^b(\Shv_1^\star)$ denote the full subcategory of $K^b(\Shv_1^\star)$ consisting
of complexes that are acyclic as complexes of $1$--motivic sheaves.
 The natural functor \[K^b(\Shv_1^\star)/N^b(\Shv_1^\star)\to D^b(\Shv_1^\fppf)\]
 is an equivalence of categories.
\end{lemma}
\proof
This follows from Lemma~\ref{lem.shv10} and \cite{KS}, Lemma 13.2.2.
\qed

Similarly, with $\Shv_1^\fr$ in place of $\Shv_1^\star$, but with no
conditions on $k$ one has
\begin{lemma}\label{lem.der1bis}
Let $N^b(\Shv_1^\fr)$ denote the full
subcategory of $K^b(\Shv_1^\fr)$ consisting of complexes that are
acyclic as complexes of $1$--motivic sheaves.
 The natural functor \[K^b(\Shv_1^\fr)/N^b(\Shv_1^\fr)\to D^b(\Shv_1^\fppf)\]
 is an equivalence of categories.
\end{lemma}

 \begin{lemma}\label{lem.der2}
 Assume condition $(*)$.
 Let $N^b(\M^\star)$ be the full subcategory of
$K^b(\M^\star)$ consisting of complexes that are acyclic as complexes of
$1$--motives with torsion.
The natural functor \[K^b(\M^\star)/N^b(\M^\star)\to D^b(\tM)\] is an equivalence of
 categories.
\end{lemma}
\proof
By Lemma~\ref{lem.mw} the ``dual" conditions required in \cite{KS}, 13.2.2 ii), are
satisfied. One checks that the ``dual" statement of \cite{KS}, 13.2.1, holds
as well, and hence one can apply \cite{KS}, 10.2.7 ii).
\qed

Similarly, with $\M$ in place of $\M^\star$, but with no conditions
on $k$ one has

\begin{lemma}\label{lem.der3}
Let $N^b(\M)$ denote the full subcategory of $K^b(\M)$ consisting of
complexes that are acyclic as complexes of $1$--motives with
torsion. The natural functor \[K^b(\M)/N^b(\M)\to D^b(\tM)\] is an
equivalence of categories.
\end{lemma}

There remains the verification that the functor ${\rm H}_0$
preserves the exactness structures.

\begin{lemma}\label{lem.acy}
 Let $M^\bullet$ be a complex in $K^b(\M^\star)$.
Then $M^\bullet \in N^b(\M^\star)$ if and only if ${\mathrm
H}_0(M^\bullet)\in N^b(\Shv_1^\star)$. In particular ${\mathrm
H}_0$ induces an equivalence of categories
\[ K^b(\M^\star)/N^b(\M^\star)\to K^b(\Shv_1^\star)/N^b(\Shv_1^\star).\]
\end{lemma}
\proof
Let $M^\bullet\colon (\cdots\to M^ i\stackrel{d^i}{\to} M^{i+1}\to\cdots)$ be a complex in
$K^b(\M^\star)$.
Observe that $\ker (d^i)\in \M^\star$ by Lemma~\ref{lem.mw}.

If $M^\bullet$ is acyclic in $K^b(\tM)$, up to isomorphisms, $\coker (d^i)$ is in $\M^\star$,
where $\coker (d^i)$ denotes here the cokernel of $d^i$ in $\tM$.
Hence, in order to prove that $H_0(M^\bullet)$ is an acyclic complex in $\Shv_1^\fppf$
it is sufficient to check the case
of a short exact sequence $M^\bullet$.
By Theorem~\ref{thm.ac} (ii) the sequence
$M^\bullet$ is represented up to q.i.\ by a sequence
 $0\to {\tilde M}^0\by{d^0}{\tilde M}^1\by{d^1} {\tilde M}^2\to 0$ of effective
$1$--motives ${\tilde M}^i=[\tilde u_i\colon \tilde L_i\to \tilde G_i]$
that is exact as sequence of complexes; moreover $\ker(\tilde u_i)=0$ by Remark~\ref{rem.ker}.
 Now, ${\mathrm H}_0(M^\bullet)$ is the sequence $0\to \coker(\tilde u_0)\to \coker(\tilde u_1)\to
\coker(\tilde u_2)\to 0$ and this is exact because of the usual ker-coker
sequence.

Let $M^\bullet\in C^{[r,s]}(\M^\star)$ and suppose that $H_0(M^\bullet)$ is
 acyclic in $K^b(\Shv_1^\fppf)$. We prove that $M^\bullet$ is acyclic
 by induction on $s-r$. Suppose first that $s-r=2$, \ie
 \[M^\bullet\colon
0\to {M}^r\by{d^{r}}{M}^{r+1}\by{d^{r+1}} {M}^s\to 0.\]
 By direct computations one sees that $d^r$ is
 a strict monomorphism, $d^{r+1}$ is an epimorphism and the cokernel of $d^r$ in $\tM^\eff$
 is q.i.\ to $M^s$. Hence $M^\bullet$ is exact in $\tM$.

Suppose now that the result is true for $s-r=n\geq 2$ and assume
$s-r=n+1$. Put
\[K^{i}=\ker({\mathrm H}_0(d^i))=\coker({\mathrm H}_0(d^{i-2}));\]
it is a sheaf in $ \Shv_1^\star$ because the category $
\Shv_1^\star$ is closed under cokernels. Let $\tilde M^i$ be the
unique (up to isomorphisms) $1$--motive in $\M^\star$ so that
${\mathrm H}_0(\tilde M^i)=K^i$. The morphism $d^{i-1}$ factors
through $\tilde M^i$ by Proposition \ref{pro.motsh} (ii). It follows
from the case $s-r=2$ that the complexes \[0\to \tilde M^{s-2}\to
M^{s-2}\to \tilde M^{s-1}\to 0, \quad 0\to \tilde M^{s-1}\to
M^{s-1}\to M^s\to 0
\] are exact, in particular $M^\bullet$ is exact at $M^s$ and $M^{s-1}$. To conclude 
one applies the induction hypothesis to the complex
$0\to M^r\to \dots\to M^{s-2}\to \tilde M^{s-1}\to 0$.
\qed

Observe that a monomorphism $\varphi$ in $\M^\star$ may not produce
a monomorphism in $\Shv_1^\star$ if the cokernel of $\varphi$ in
$\tM$ is not isomorphic to an object in $\M^\star$. As a
counterexample consider the $p$-multiplication on $\G_m$.

In view of the foregoing lemmas our main result generalizing \cite{BVK}, 1.6.1
and 3.9.2 now follows immediately:

\begin{thm}\label{thm.main} Set $D^b(\M)\df K^b(\M)/N^b(\M)$. If $k$ has characteristic zero or if it is
transcendental over its prime field, we have canonical equivalences
of categories
\[D^b(\M)\cong D^b(\tM)\cong D^b(\Shv_1^\fppf).\]
\end{thm}

\begin{remark}\label{rem.restriction}
We have already seen
 that for finite fields, or more generally algebraic extensions of finite
fields, the categories $\M^\star$ and $\Shv_1^\star$ are both equivalent
to the category of semiabelian varieties. They then have ``not
enough objects'' to encode all information needed to describe the
bounded derived categories of $\tM$ and $\Shv_1^\fppf$. By Lemmas
\ref{lem.der1bis} \& \ref{lem.der3} the categories $\M$ and
$\Shv_1^\fr$ would be good candidates over any field; unfortunately
we do not see any equivalence between them.
\end{remark}

\begin{remark} For $X$ a smooth projective $k$-variety
the sheaf ${\mathrm {Pic}}_{X/k}$ is clearly $1$--motivic since it is
representable by a group scheme of finite type whose reduced
identity component is an abelian variety. In \cite{BVK}, 3.4.1, the
authors prove that the relative Picard functor is $1$--motivic for
the \'etale topology as soon as $X$ is smooth over $k$, \ie with the
notations in (\ref{pic.sh}) $\pi_*{\mathrm {Pic}}_{X/k}$ is a sheaf
in $\Shv_1^\et$. Unfortunately the proof in \cite{BVK}, 3.4.1, does not
work in the fppf context and it is not clear at the present if a
similar result holds for the fppf topology.
\end{remark}

\section{Further results on $1$--motivic sheaves.}\label{sec.result}

In this section 
we describe the relation between our definition of
$1$--motivic sheaves and the definition given in
\cite{BVK}. Furthermore we present
 an alternative definition of $1$--motivic sheaves
as cokernels of morphisms of group schemes \mbox{$F_1\to F_0$} where $F_1$
is an object in $\mathcal {CE}$, \ie product of a discrete group by
a finite connected commutative group scheme, and $F_0$ is a
commutative $k$-group scheme extension of an object of $\mathcal
{CE}$ by a semiabelian variety.

\subsection*{Comparison of topologies}

Let $\Shv_1'$ be the full subcategory of $\Shv_1^\fppf$ whose objects are
sheaves $\mathcal F$ as in (\ref{eq.F}) with $L, E$ discrete. It coincides with
$\Shv_1^\fppf$ in characteristic $0$. 
Denote by {\bf Sm}$/k$ the category of smooth separated $k$-schemes.
Let $\Shv_1^\et$ be the category of \'etale sheaves on ${\bf Sm}/k$
 that fit in a sequence like (\ref{eq.F}) with $G$ semiabelian,
$L$, $E$ discrete.
 The definition of the category of $1$--motivic sheaves $\Shv_1$
in \cite{BVK} is, with our notations,
\[\Shv_1\df \Shv_1^\et[1/p] .\]
Denote by $\Shv_1^{\et,\star}$ the full subcategory of $\Shv_1^\et$
consisting of those sheaves with $E=0$. Let $\pi\colon ({\rm
Sch}/k)_{\fppf}\to ({\bf Sm}/k)_\et $ be the usual morphism of
sites and consider the following diagram
\begin{eqnarray}\label{pic.sh}
\xymatrix{
\Shv_1^\star \ar[r] &\Shv_1'\ar[r]& \Shv_1^\fppf \ar[r]& \Sh(({\rm Sch}/k)_{\fppf}) \\
\Shv_1^{\et, \star} \ar[u]^{\pi^*}\ar[r] &\Shv_1^\et\ar[u]^{\pi^*} \ar[rr]&
 & \Sh(({\bf Sm}/k)_{\et})\ar[u]^{\pi^*}
 }
\end{eqnarray}

\begin{lemma}
The functors $\pi^*\colon \Shv_1^\et\to \Shv_1'$ and 
$\pi^*\colon \Shv_1^{\et, \star}\to \Shv_1^\star$ are equivalences of
categories with quasi-inverse $\pi_*$.
\end{lemma}
\proof
(See also \cite{BVK}, 3.3.2).
First of all observe that for $X$ a smooth $k$-scheme one has $\pi^*X=X$
(proof as in \cite{EtC}, p.~69). Furthermore, $\pi_*\pi^*X=X$.
Indeed, let $U$ be
a smooth $k$-scheme. Then $\Gamma(U,\pi_*\pi^*X)=\Gamma(U,\pi^*X)$ by
definition of $\pi_*$; we have just seen that the latter group
equals $\Gamma(U,X)$ and this does not change when working with the
\'etale or the flat topology.

Moreover $R^i\pi_*X=0$, $i>0$, for $X$ a smooth $k$--group scheme.
Indeed $R^i\pi_*X=0$ if for all smooth $k$--schemes $S$ one has
${\rm H}^i(S_\et,X)={\rm H}^i(S_\fppf,X)$, and this last fact follows from
\cite{GRO}, 11.7. In particular, for $\mathcal F$ in $\Shv_1'$, one has $R^i\pi_*(G/L)=0=R^i\pi_*({\mathcal F})$,
 $i>0$.

We start by showing that $\pi^*\pi_*$ is naturally isomorphic to the identity on $\Shv_1'$. Let
$\mathcal F$ be a $1$--motivic sheaf in $\Shv_1'$. By
Proposition~\ref{pro.van}, $\mathcal F=F_0/F_1$ with $F_0,F_1$ smooth group
schemes over $k$, $F_1$ discrete. Since $R^1\pi_*F_1=0$, we have that
$\pi_*{\mathcal F}=\pi_*F_0/\pi_*F_1=F_0/F_1$ is an object of $\Shv_1^\et$
(\cf~Remark~\ref{rem.et}).
Now, $\pi^*\pi_*{\mathcal F}=\pi^*\pi_*F_0/\pi^*\pi_*F_1$ by the right
exactness of $\pi^*$. In particular, $\pi^*\pi_*{\mathcal F} ={\mathcal F}$.

 We now show that $\pi_*\pi^*$ is naturally isomorphic to the identity on
 $\Shv_1^\et$. Let $\mathcal F$ be an object in $\Shv_1^\et$.
By Remark~\ref{rem.et}, $\mathcal F=F_0/F_1$ with $F_0,F_1$ smooth $k$--group
schemes.
Now, since $\pi^*$ is right exact,
$\pi^*{\mathcal F}=\pi^*F_0/\pi^*F_1=F_0/F_1$ is an object of $\Shv_1'$ (\cf~Proposition~\ref{pro.van}). Furthermore,
since $R^1\pi_*F_1=0$, ${\mathcal F}=F_0/F_1=$$\pi_*\pi^*F_0/\pi_*\pi^*F_1=$$\pi_*\pi^*{\mathcal F}$.
\qed

\begin{propose}\label{pro.voe}\begin{itemize}
\item[(i)] The category $\Shv_1^\et$ (respectively $\Shv_1^{\et,\star}$) is equivalent to the full subcategory $\Shv_1'$ (respectively
$\Shv_1^\star$) of $\Shv_1^\fppf$. The subcategory $\Shv_1'$
is cogenerating and closed under
cokernels.
\item[(ii)] Denote by $N^b(\Shv_1')$ the bounded complexes of
objects in $\Shv_1'$ that are acyclic as complexes of $1$--motivic sheaves.
 The natural functor \[K^b(\Shv_1')/N^b(\Shv_1')\to D^b(\Shv_1^\fppf)\]
 is an equivalence of categories.
\item[(iii)] If ${\rm char} k=0$ or if $k$ is transcendental over its prime field, we have an equivalence of categories.
\[ K^b(\Shv_1^{\et})/N^b(\Shv_1^{\et })\to D^b(\Shv_1^\fppf)\cong D^b(\M). \]
\end{itemize}
\end{propose}
\proof The equivalence results in (i) were proved above. Since
$\Shv_1'$ contains $\Shv_1^\star$, it is cogenerating. 
Hence (ii) can be proved in the same way as Lemma~\ref{lem.der1}. 
For (iii), we have the following equivalences
\begin{multline*}K^b(\Shv_1^{\et})/N^b(\Shv_1^{\et })\cong
 K^b(\Shv_1^{\et,\star})/N^b(\Shv_1^{\et,\star })\cong
 K^b(\Shv_1^\star)/N^b(\Shv_1^\star)\\
 \cong D^b(\Shv_1^\fppf)
 \cong D^b(\M),
\end{multline*}
 where the first equivalence is due to \cite{KS}, 10.2.7, and the second
 equivalence is due to the fact that both $\pi^*$ on $\Shv_1^{\et,\star}$ and $\pi_*$ on
$\Shv_1^\star$ are exact functors. The remaining equivalences were
proved in Lemma~\ref{lem.der1} and Theorem~\ref{thm.main}.
\qed

\begin{remark} Observe that $D^b(\Shv_1^\fppf[1/p])$ is equivalent to the thick subcategory
$d_{\leq 1}{\rm DM}^\eff_{\rm{gm},\et}$ of Voevodsky's triangulated category of
motives generated by motives of smooth curves.
Indeed the
inclusion functor $\Shv_1'\to \Shv_1^\fppf$ provides an equivalence
of categories $\Shv_1'[1/p]\to \Shv_1^\fppf[1/p]$ because any sheaf
in $\Shv_1^\fppf$ associated to a finite group scheme of order a
power of $p$ is isomorphic to $0$ in $\Shv_1^\fppf[1/p]$; in
particular, we are killing the connected component of $E$ in
(\ref{eq.F}). To conclude one recalls that $\Shv_1'[1/p]$ is equivalent to $
\Shv_1^\et[1/p]$ by Proposition \ref{pro.voe} (i) and that
 $D^b(\Shv_1^\et[1/p])$ is equivalent to $d_{\leq 1}{\rm
 DM}^\eff_{\rm{gm},\et}$ by \cite{BVK}, 3.9.2. 
Observe that the latter category is $\Z[1/p]$-linear by \cite{Voe}, 3.3.3.
\end{remark}

\subsection*{Alternative definition of $1$--motivic sheaves}

We now see that any $1$--motivic sheaf admits a presentation as
cokernel of a morphism $F_1\to F_0$ of group schemes where $F_1$ is
discrete and the reduced subgroup of the identity component of $F_0$
is semiabelian. This fact was used to prove Lemma \ref{lem.shv10} and
Proposition \ref{pro.voe}.
We start by showing some vanishing results for $2$-fold extensions.

\begin{lemma}\label{lem.extNL}
Let $N$ be a finite connected commutative group scheme over $k$ and $L$ a
discrete group. Then $\Ext^2(N,L)=0$
\end{lemma}
\proof Let $n$ be the order of $N$ and $\Ext^1(N,L/nL)\to \Ext^2(N,L)$
the boundary map deduced from the exact sequence $0\to L\to L\to
L/nL\to 0$. This homomorphism is surjective and the first group is trivial over a perfect field.
\qed
\begin{lemma}\label{lem.ext2EL}
Let $E,L$ be discrete groups. Any $2$-fold extension as in \eqref{eq.F}
becomes trivial after pull-back along a suitable epimorphism $\tilde
E\to E$ of discrete groups.
\end{lemma}
\proof It is sufficient to prove that any extension 
$\eta\colon 0\to G/L\to {\mathcal F}\to E\to 0$ splits after pull-back along a
suitable epimorphism of discrete groups $\tilde E\to E$ 
(\cf~\cite{BVK}, 3.7.5, Step~1).

First of all observe that $E$ becomes constant after base change 
to a suitable finite Galois extension $k'$ of $k$, 
 and that if we write $f\colon \Spec(k')\to\Spec(k)$ for the
 corresponding morphism of schemes, the base change functor $f^*$ is exact, 
\ie $f^*\eta$ provides a $1$--motivic sheaf over $k'$.
 Up to enlarging $k'$, we may assume that the extension
 $\eta$ splits over $k'$:
if $E_{k'}$ is torsion free, this follows from the fact that any
element in $\mathrm{ Ext}^1(\Z, G)={\mathrm H}_\fppf^1(k',G)=
{\mathrm H}_\et^1(k',G)$, becomes trivial
 over a suitable finite Galois extension of $k'$ and similarly for
${\mathrm H}_\et^2(k',L) $; 
if $E_{k'}$ is $n$-torsion, one
 uses the fact that any element in $\mathrm{ Ext}^1(\Z/n\Z, G)$, as well as
 any element in $\mathrm{ Ext}^2(\Z/n\Z, L)$, becomes trivial over a suitable
finite Galois extension of $k'$.
Let $g\colon E_{k'}=f^*E\to
 f^*\mathcal F$ be a splitting of $\eta$ over $k'$.
Let $$ {\mathcal N}\colon \prod_{\sigma\in Gal(k'/k)} E_{k'}^\sigma =f^*f_*E_{k'}\to
 E_{k'}, \quad
 (a_\sigma)_\sigma\to \sum_\sigma a_\sigma.$$ 
This homomorphism descends to the trace map (\cf~\cite{EtC}, V, 1.12) 
 $f_*E_{k'}\to E$ which is an epimorphism of discrete groups. Furthermore,
$g\circ \mathcal N$ descends to a trivialization of the pull-back of $\eta$ along the trace map.
 \qed

The following proposition provides an alternative definition of
$1$--motivic sheaves.

\begin{propose}\label{pro.van}
An fppf sheaf $\mathcal F$ on $k$ is
\emph{$1$--motivic} if and only if
\begin{eqnarray}\label{eq.motsh}
\mathcal F=\coker(F_1\by{u} F_0)
\end{eqnarray}
where $F_1$ is a discrete group scheme, $F_0$ is an extension of an object
in $\mathcal {CE}$ by a semiabelian group scheme $G$ and $u$ is a monomorphism.

Denote by $L$ the pull-back
of $F_1$ to $G$; we then have a diagram
\begin{eqnarray}\label{dia.def}
\xymatrix{
0\ar[r]& L\ar[r]\ar@{^{(}->}[d]& G\ar[r]\ar[dr]^b\ar@{^{(}->}[d]& \mathcal
F^\star\ar[r]\ar@{^{(}->}[d]& 0\\
 0\ar[r] & F_1 \ar@{->>}[d] \ar[r]^u & F_0 \ar[r]\ar@{->>}[d] & {\mathcal F}\ar[r]\ar@{->>}[d]& 0\\
0\ar[r]& F_1/L \ar[r] & F_0/G \ar[r]& E\ar[r] & 0
 }\end{eqnarray}
For a $1$--motivic sheaf $\mathcal F$ one can always find a presentation
$F_1\to F_0$ as above with $F_1$ \'etale and torsion free.
If $k$ is algebraically closed or $E$ is connected,
there exists a diagram as above with $F_1/L=0$.
\end{propose}
\proof
 The if part follows from diagram~\eqref{dia.def}.
For the converse, suppose $b\colon G\to \mathcal F$ normalized.
Let $\eta$ be the $2$-fold extension in \eqref{eq.F}.
If $\eta$ is isomorphic to the trivial extension, then $\mathcal F$ is
(isomorphic to)
 the push-out along $G\to G/L=:{\mathcal F}^\star$ of
an extension $F_0$ of $E$ by $G$ and we get a diagram as above with
$F_1/L=0$. By Lemma~\ref{lem.extNL} this is the case if $E$ is
finite connected. There remains only to check the case $E=E^\et$ \'etale. By
Lemma~\ref{lem.ext2EL}, there exists an epimorphism $\tilde E\to E$ such that
the class of the pull-back $0\to L\to G\to
\tilde {\mathcal F}\to \tilde E\to 0$ of $\eta$ along $\tilde E\to
E$ is trivial and then, as explained above, $\tilde {\mathcal F}=
\coker (L\to F_0)$ with $ F_0$ extension of $\tilde E$
by $G$. Furthermore the composition $ F_0\to \tilde {\mathcal
F}\to \mathcal F$ remains an epimorphism whose kernel $F_1$ is
extension of $\ker(\tilde E\to E)$ by $L$. \qed

\begin{remark}\label{rem.et}
The previous Proposition also holds for the categories $\Shv_1^\et$ and
$\Shv_1^\et[1/p]$. Hence an \'etale sheaf $\mathcal F$ on ${\bf Sm}/k$ is
$1$--motivic in the sense of \cite{BVK} if and only if $\mathcal F$
 fits into a diagram as above with $F_1$ discrete and 
$F_0$ a smooth commutative $k$-group with semiabelian identity component. 
Indeed in \cite{BVK}, 3.7.5, it is shown that one can
always find such a presentation with $F_0$ a split extension of its
component group by its identity component. Furthermore, in
\cite{BVK}, 1.3.8, it is shown that the \'etale sheaves associated to
the group schemes $F_0, F_1$ in \eqref{dia.def} are objects of the
category ${\rm HI}_\et$, \ie homotopy invariant \'etale sheaves with
transfers over $k$. One can then construct a functor $\rho\colon
\Shv_1^\et \to {\rm HI}_\et$ mapping $\mathcal F$ to $\coker(F_1\to
F_0)$. As a sheaf, this cokernel is still $\mathcal F$ and by
\cite{BVK}, 3.8.2, its transfer structure does not depend on the
presentation. Inverting $p$-multiplications, the functor $\rho$
induces the full embedding $\Shv_1:=\Shv_1^\et[1/p] \to {\rm
HI}_\et[1/p]$ studied in \cite{BVK}, 3.8.1.
\end{remark}

{\bf Acknowledgments:}
The author would like to thank L. Barbieri-Viale for
pointing her attention to this subject and for helpful discussions.

 \end{document}